\documentclass[11pt]{article}
\usepackage{amssymb,epsfig}

\newtheorem{theorem}{Theorem}[section]
\newtheorem{proposition}[theorem]{Proposition}

\newtheorem{example}[theorem]{Example}

\newenvironment{proof}{{\noindent \sc Proof. }}{\hfill $\Qed$}

\newcommand{\la}{\langle}
\newcommand{\ra}{\rangle}

\newcommand{\Qed}{\rule{2.5mm}{3mm}}

\newcommand{\Aut}{\hbox{{\rm Aut}}\,}

\newcommand{\ZZ}{\mathbb{Z}}

\newcommand{\C}{{\cal{C}}}

\newcommand{\HH}{{\cal{H}}}

\begin{document}

%% title page %%%%%%%%%%%%%%%%%%%%%%%%%%%%%%%%%%%%%%%%%%%%%%%%%%%%%

\begin{center}
{\bf\large HAMILTONICITY OF CUBIC CAYLEY GRAPHS}
\end{center}

\bigskip\noindent
\begin{center}
\begin{tabular}{cc}
{\sc Henry Glover}
& {\sc Dragan Maru\v si\v c} \footnotemark \\
{\small Department of Mathematics}
& {\small IMFM, University of Ljubljana} \\
{\small Ohio State University}
& {\small and University of Primorska, Koper} \\
{\small U.S.A.}
& {\small Slovenia} \\
{\small {\tt glover@math.ohio-state.edu}}
& {\small {\tt dragan.marusic@guest.arnes.si}} \\
\end{tabular}
\end{center}

\addtocounter{footnote}{0}
%\footnotetext{Supported in part by ``Ministrstvo za \v solstvo,
%znanost in \v sport Republike Slovenije'', research project Z1-3124.}
%\addtocounter{footnote}{1}
\footnotetext{Supported in part by ``Ministrstvo za visoko \v solstvo,
znanost in tehnologijo Republike Slovenije'', research program P1-0285.}

\begin{abstract}

 Following a problem posed by Lov\'asz in 1969,
 it is believed that every connected vertex-transitive graph has a Hamilton path.
 This is shown here to be true for cubic Cayley graphs arising from
 groups  having  a $(2,s,3)$-presentation, that is, for groups
 $G=\la a,b| a^2=1, b^s=1, (ab)^3=1, etc. \ra$ generated by an involution
 $a$ and an element $b$ of order $s\geq3$
 such that their product $ab$ has order $3$. More precisely, it is
 shown that the Cayley graph $X=Cay(G,\{a,b,b^{-1}\})$
 has a Hamilton cycle when  $|G|$ (and thus $s$) is congruent to $2$ modulo $4$,
 and has a long cycle missing only two vertices (and thus
 necessarily a Hamilton path) when $|G|$ is congruent to $0$ modulo $4$.
 \end{abstract}

%
%
%%%%%%%%%%%%%%%%%%%%%%%%%%%%%%% section 1 %%%%%%%%%%%%%%%%%%%%%%%%%%%%%%%%%%%%%
%
%

\section{Introductory remarks}
\label{sec:intro}
\noindent

\medskip

In 1969, Lov\'asz \cite{LL70} asked whether 
every connected vertex-transitive graph has a Hamilton path,
thus tying together, through this special case of the Traveling Salesman Problem,
two seemingly unrelated concepts: traversability and symmetry of graphs.
Lov\'asz problem is, somewhat misleadingly, usually 
referred to as the Lov\'asz conjecture,
presumably in view of the fact that, after all these years,
a connected vertex-transitive graph without a Hamilton path
is yet to be produced. Moreover, only four connected
vertex-transitive graphs (having at least three vertices)
not possessing a Hamilton cycle
are known to exist:
the Petersen graph, the Coxeter graph, and the two graphs obtained
from them by replacing each vertex with a triangle.
All of these are cubic graphs,
suggesting perhaps that no attempt to resolve the above problem
can bypass a thorough analysis of cubic vertex-transitive graphs.
Besides, none of these four graphs is a Cayley graph.
This has led to a folklore conjecture that every Cayley graph is hamiltonian.

This problem has spurred quite a bit of interest in the mathematical community.
In spite of a large number of articles directly and indirectly related
to this subject
(see \cite{BA79,BA83,AP82,ADP85,BA89,AZ89,BA89a,ALW90,AM96,AQ01,YC98,DGMW98,ED83,GR85,GY96,
KW85,DP82,DP83,DM83,DM85,DM87,DM88,DW82,DW85,DW86} for some
of the relevant references), not much progress has been made
with regards to either of the two conjectures.

For example, most of the results proved thus far in the case of Cayley graphs
depend on various restrictions made either on the class of groups dealt with or
on the generating sets of Cayley graphs.
For example, one may easily see that Cayley graphs of abelian groups
have a Hamilton cycle. Also, following a series of articles \cite{ED83,KW85,DM83} it is now
known that every Cayley graph of a group with a cyclic commutator subgroup of prime power order,
is hamiltonian. This result has later been generalized
to connected vertex-transitive graphs
whose automorphism group contain a transitive subgroup whose
commutator subgroup is cyclic of prime-power order,
with the Petersen graph being the only counterexample \cite{DGMW98}.
And finally, perhaps the biggest achievement on the subject
is a result of Witte (now Morris) which says that
a Cayley (di)graph of any $p$-group has a Hamilton cycle \cite{DW86}.
(For further  results not explicitly mentioned or referred to here
see the survey paper \cite{CG96}.)

In this  article we consider the hamiltonicity problem for
cubic Cayley graphs arising from groups
having  a $(2,s,3)$-presentation, that is, for groups
$G=\la a,b| a^2=1, b^s=1, (ab)^3=1, etc. \ra$ generated by an involution
$a$ and an element $b$ of order $s\geq3$
such that their product $ab$ has order $3$.
More precisely, the following is the main result of this article.

\begin{theorem}
\label{the:main}
Let $s\geq3$ be an integer and let $G=\la a,b| a^2=1, b^s=1, (ab)^3=1, etc. \ra$
be a group with a $(2,s,3)$-presentation. Then the
Cayley graph $X=Cay(G,\{a,b,b^{-1}\})$ has a
Hamilton cycle when $|G|$ (and thus also $s$) is congruent to $2$ modulo $4$,
and has a cycle of length $|G|-2$, and thus
necessarily a Hamilton path,  when $|G|$ is congruent to $0$ modulo $4$.
\end{theorem}

Let us comment that the class of groups considered in
Theorem~\ref{the:main} is by no means restrictive.
First, by \cite{MSW94}, \cite{LM99}, \cite{AS98} and \cite{AW89} every finite
nonabelian simple
group except the Suzuki groups, a thin family of $PSp_n(q)$ and
a thin family of $PSU_n(q)$ groups, $M_{11}$, $M_{22}$, $M_{23}$, $McL$
and at most finitely many other non-sporadic finite simple groups have
a $(2,s,3)$-presentation. Also, methods similar to those in this article
have been used in \cite{CS94,CSW89} to find Hamilton cycles in certain
Cayley graphs.
And second, if $X$ is a cubic arc-transitive
graph and $G \leq \Aut X$ acts $1$-regularly on $X$,
then it is easily seen that $G$ has a $(2,s,3)$-presentation for some $s$.
Namely, the ordered pair $(X,G)$ gives rise to a unique orbit 
of those undirected cycles in $X$ which have the property that 
each of these cycles is rotated by some automorphism in $G$
(that is, the so called {\em consistent cycles} in the terminology of Biggs \cite{NB78}).
These cycles give rise to the faces of the corresponding
(orientably) regular map associated with $X$, and their length  
is then precisely our parameter $s$ in the $(2,s,3)$-presentation of $G$.
Going backwards, the well defined correspondence
between these groups (or rather their Cayley graphs)
and the class of all those cubic arc-transitive
graphs which admit a subgroup acting regularly on the arcs is,
geometrically, best seen via the concept of the hexagon graphs, explained
in the subsequent section.
(However, this correspondence is not 1-1, for a cubic arc-transitive graph
may possess nonisomorphic $1$-regular subgroups.)

The article is organized as follows.
In Section~\ref{sec:ex} we describe our method for constructing Hamilton
cycles and paths in cubic Cayley graphs of groups having a $(2,s,3)$-presentation
by analyzing six examples of such graphs.
They are associated with, respectively,
the groups $\ZZ_6$ and $S_3 \times \ZZ_3$ with a $(2,6,3)$-presentation,
the group $S_4$  with a $(2,4,3)$-presentation,
the group $Q_8 \rtimes S_3$  with a $(2,8,3)$-presentation,
the group $A_4$ with a $(2,3,3)$-presentation
and the group $A_5$ having a $(2,5,3)$-presentation.
In Section~\ref{sec:cyclic} we introduce the graph-theoretic concepts of cyclic stability
and cyclic connectivity. In particular, we discuss an old theorem
of Payan and Sakarovitch \cite{PS75} which gives the exact size
of a maximum cyclically stable set in a cyclically $4$-connected
cubic graph (Proposition~\ref{pro:pasa}),
a result that proves to be of crucial importance
for the purpose of this article.
Using a result of Nedela and \v Skoviera \cite{NS95}
on the cyclic connectivity in cubic vertex-transitive graphs
(Proposition~\ref{pro:nesko}),
together with an analysis of cubic arc-transitive graphs
of girth at most $5$ (Proposition~\ref{pro:girth}),
we then obtain a slight refinement
of the above mentioned result of Payan and Sakarovitch
(Proposition~\ref{pro:modi}), thus laying the groundwork for the proof
of Theorem~\ref{the:main} which is carried out in Section~\ref{sec:proof}.

%
%
%%%%%%%%%%%%%%%%%%%%%%%%%%%%%%% section 2 %%%%%%%%%%%%%%%%%%%%%%%%%%%%%%%%%%%%%
%
%

\section{The method of proof illustrated}
\label{sec:ex}
\noindent

\medskip

In this section we give examples illustrating our method of proof of Theorem~\ref{the:main}.
In particular, each Cayley graph we study has a canonical Cayley map given by
an embedding of the Cayley graph $X=Cay(G,\{a,b,b^{-1}\})$
of the (2,s,3)-presentation of  a group $G=\la a,b| a^2=1, b^s=1, (ab)^3=1, etc. \ra$
in the closed
orientable surface of genus $1+(s-6)|G|/12s$ with faces $|G|/s$ disjoint $s$-gons
and $|G|/3$ hexagons. This map is given by using the same rotation of
the $b$, $a$, $b^{-1}$ edges at every vertex and results in one $s$-gon and
two hexagons adjacent to each vertex. In each case we give a tree of either
$(|G|-2)/4$ hexagons if $|G|\equiv 2(mod\,4)$ and $(|G|-4)/4$ hexagons if $|G|\equiv 0(mod\,4)$.
This tree of hexagons necessarily contains, respectively, all or all but two of the
vertices of the Cayley graph and as a subspace of the Cayley map is a
topological disk. The boundary of this topological disk is a (simple) cycle
passing through, respectively, all or all but two vertices of the Cayley graph.
We give two examples in the case $|G|\equiv 2(mod\,4)$ and four examples in the case
$|G|\equiv 0(mod\,4)$.
In each example we show the tree of hexagons in the Cayley map,
in the first case giving rise
to a Hamilton cycle of the graph and in the second case giving rise to
a long cycle missing only two vertices.
Finally we do show a Hamilton cycle in the Cayley graph in the
case $|G|\equiv 0(mod\,4)$ when the tree of $(|G|-4)/4$ hexagons does not give a
Hamilton cycle in the Cayley graph. We do this by exhibiting a {\em Hamilton tree}
of faces in the Cayley map (a tree of faces such that each vertex of the Cayley
graph lies in the boundary of at least one of these faces) by
using an appropriate number of $s$-gons.
In Examples~\ref{ex:z6} and \ref{ex:s3z3}
we have  $|G|\equiv 2(mod\,4)$ and $s\equiv 2(mod\,4)$,
in Examples~\ref{ex:s4} and \ref{ex:q8s3}
we have  $|G|\equiv 0(mod\,4)$ and $s\equiv 0(mod\,4)$,
and in Examples~\ref{ex:a4} and \ref{ex:a5}
we have $s\equiv 1(mod\,2)$ and thus $|G|\equiv 0(mod\,4)$.
%(*** as was pointed to us by Marston Conder.)

Note that the above construction has a direct translation into
a more graph-theoretic language by associating
with the Cayley graph $X=Cay(G,\{a,b,b^{-1}\})$ of $G$ the so called
{\em hexagon graph} $Hex(X)$ whose vertex set
consists of all the hexagons in $X$ arising from the relation
$(ab)^3$, with  two hexagons adjacent in $Hex(X)$
if they share an edge in $X$.
It may be easily seen that $Hex(X)$ is nothing but
the so called {\em orbital graph} of the left action of $G$ on the set $\HH$
of left cosets of the subgroup $H=\la ab \ra$, arising from
the suborbit $\{aH,abaH, ababaH\}$ of length $3$.
(But note that $aH=bH$ and so $abaH=ab^2H$ and $ababaH =b^{-1}H.$)
More precisely, the graph has  vertex set
$\HH$, with adjacency defined as follows:
an arbitrary  coset $xH$ is adjacent to
precisely the three cosets
$xbH$, $xb^{-1}H$ and $xab^2H$.
Clearly, $G$ acts $1$-regularly on $Hex(X)$.
Conversely, let $X$ be a cubic arc-transitive graph $Y$
admitting a 1-regular action of a subgroup $G$ of $\Aut Y$.
Let $v \in V(Y)$ and let $h$ be a generator of $H =G_v \cong \ZZ_3$.
Then there must exist an element $a \in G$ such that
$G = \la a,h \ra$ and such that $Y$ is isomorphic to the
orbital graph of $G$ relative to the suborbit
$\{aH, haH, h^2aH\}$. Moreover, a short computation
shows that $a$ may be chosen to be an involution,
and letting $b =ah$ we get the desired presentation for $G$.
There is therefore a well defined correspondence
between these two classes of objects, as
noted in the introductory section.
However, this correspondence is not 1-1, for a cubic arc-transitive graph
may possess nonisomorphic $1$-regular subgroups.
A typical example is the Moebius-Kantor graph on $16$
vertices which admits two noisomorphic
$1$-regular subgroups, one with a $(2,8,3)$-presentation
and the other with a $(2,12,3)$-presentation.
The former and the corresponding Cayley graph
is discussed in Example~\ref{ex:q8s3} below.

The trees of hexagonal faces in the associated Cayley map of $X$
(mentioned in the first paragraph)
then correspond to vertex subsets in $Hex$
inducing trees with the property that
the complement in $V(X)$ is either an independent set
when $|G| \equiv 2\,(mod\,4)$,
or induces a subgraph with  a single edge
when $|G| \equiv 0\,(mod\,4)$.
That this approach works in general will follow
from the results given in Section~\ref{sec:cyclic}.

\begin{figure}[h]
\begin{center}
%\label{fig:z6}
\includegraphics[width=15cm]{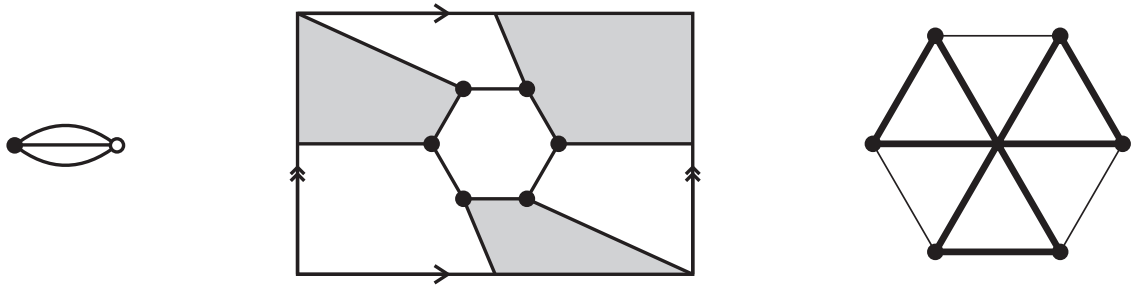}
\caption{A (trivial) Hamilton tree of faces in a toroidal Cayley
map of of $K_{3,3}$ giving rise to a Hamilton cycle, and the
associated hexagon graph.}\label{fig:z6}
\end{center}
\end{figure}

\begin{example}
\label{ex:z6}
{\rm In the middle picture of  Figure~\ref{fig:z6} we show
a trivial tree of hexagons (consisting of a single hexagon),
whose boundary is a Hamilton cycle
in the toroidal Cayley map
of $X=K_{3,3}$, the Cayley graph of the group $G = \ZZ_6$
with  a $(2,6,3)$-presentation $\la a,b \mid a^2=b^6=(ab)^3 =1, etc. \ra$,
where $a=3$ and $b=1$. The left picture shows the corresponding
hexagon graph $\Theta_2$, and the right picture shows the
corresponding Hamilton cycle in $X$.}
\end{example}

\bigskip
%\bigskip

\bigskip
%\bigskip

\begin{figure}[h]
\begin{center}
%\label{fig:s3z3}
\includegraphics[width=15cm]{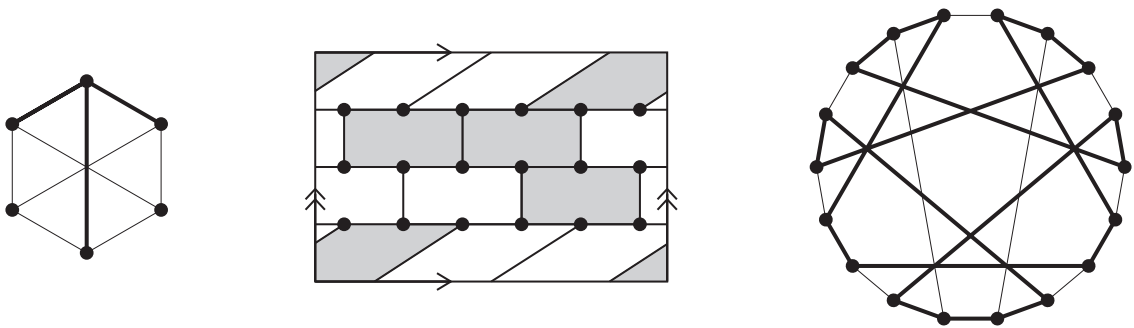}
\caption{A Hamilton tree of faces in a toroidal Cayley map of the
Pappus graph giving rise to a Hamilton cycle, and the associated
hexagon graph.}\label{fig:s3z3}
\end{center}
\end{figure}

\begin{example}
\label{ex:s3z3}
{\rm In the middle picture of  Figure~\ref{fig:s3z3} we show
a Hamilton tree of hexagons, whose boundary is a Hamilton cycle
in the toroidal Cayley map of
the Pappus graph $X$, a Cayley graph of the group $G = S_3 \times \ZZ_3$
with  a $(2,6,3)$-presentation $\la a,b \mid a^2=b^6=(ab)^3 =1, etc. \ra$,
where $a=((12),0)$ and $b=((13),1)$.
The left picture shows this same tree in the corresponding
hexagon graph $K_{3,3}$,
and the right picture shows the
corresponding Hamilton cycle in $X$.}
\end{example}

%\bigskip
\bigskip

\bigskip
\bigskip

\begin{figure}[h]
\begin{center}
%\label{fig:s4}
\includegraphics[width=15cm]{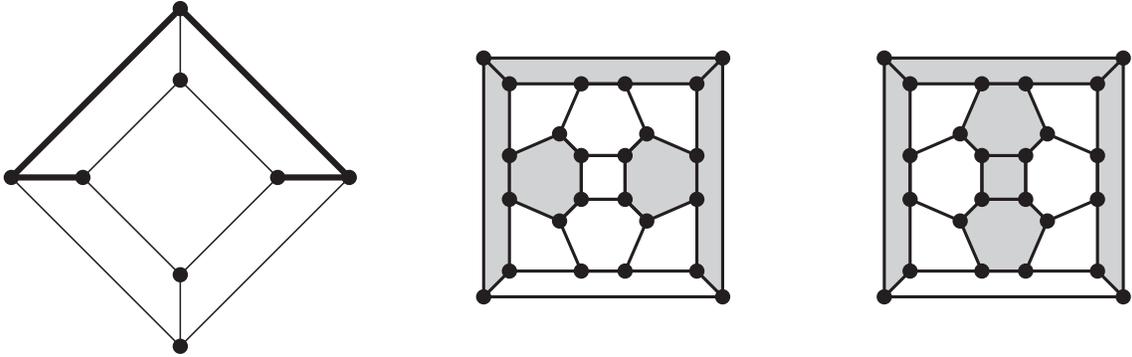}
\caption{A tree of faces in the spherical Cayley map of a Cayley
graph of $S_4$ giving rise to a cycle missing two vertices, the
associated hexagon graph, and a modified Hamilton tree of
faces.}\label{fig:s4}
\end{center}
\end{figure}

\begin{example}
\label{ex:s4}
{\rm In the middle picture of  Figure~\ref{fig:s4} we show
a tree of hexagons, whose boundary is a cycle missing only two vertices
in the spherical Cayley map of
a Cayley graph $X$ of the group $G = S_4$
with  a $(2,4,3)$-presentation $\la a,b \mid a^2=b^4=(ab)^3 =1 \ra$,
where $a=(12)$ and $b=(1234)$.
The left picture shows this same tree in the corresponding
hexagon graph $Q_3$, the cube,
and the right picture shows a modified tree of faces,
including also a square, whose boundary is a Hamilton cycle
in this map.}
\end{example}

\bigskip
\bigskip

\bigskip
\bigskip

\begin{figure}[h,t]
\begin{center}
%\label{fig:q8s3a}
\includegraphics[width=14cm]{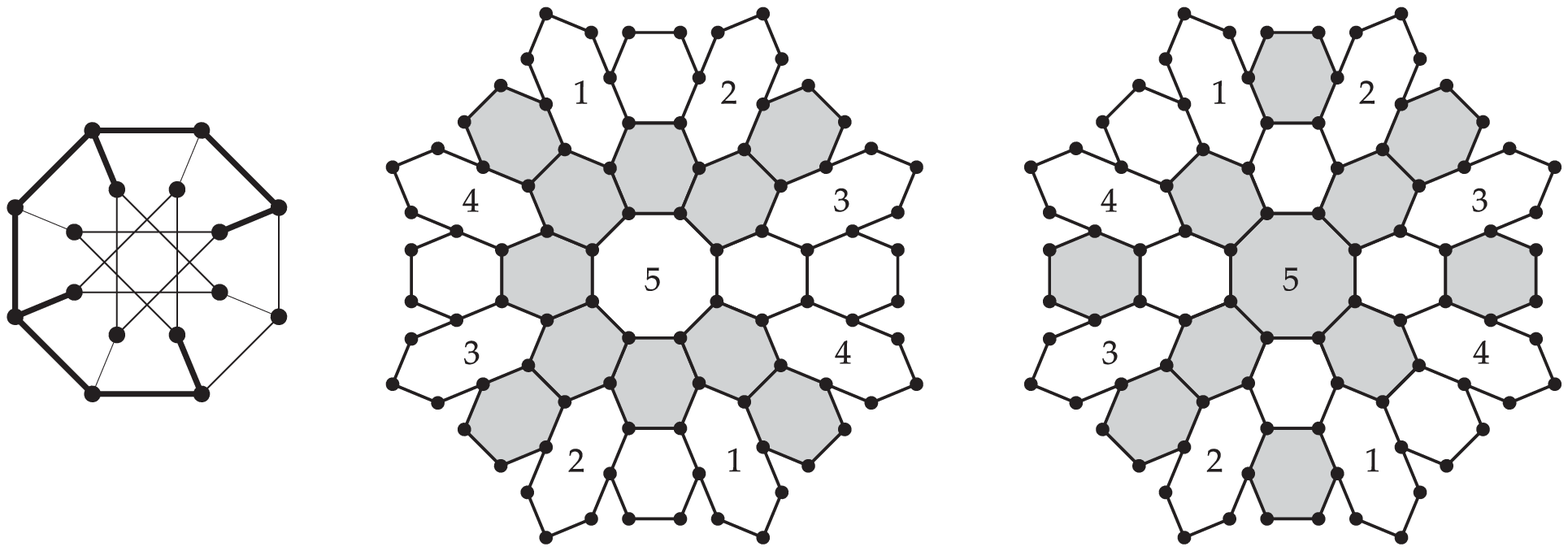}
\caption{A tree of faces in the genus $2$ Cayley map of a
Cayley graph of $Q_8 \rtimes S_3$ giving rise to a cycle missing
two vertices, the associated hexagon graph, and a modified
Hamilton tree of faces.}\label{fig:q8s3a}
\end{center}
\end{figure}

\begin{example}
\label{ex:q8s3} {\rm In the middle picture of
Figure~\ref{fig:q8s3a} we give the genus $2$ Cayley map of a
Cayley graph $X$ of the group $G = Q_8 \rtimes S_3$ with  a
$(2,8,3)$-presentation $\la a,b \mid a^2=b^8=(ab)^3 =1, etc. \ra$,
where $a = (1,(23))$ and $b=(i,(12))$. The action of the
transposition $(12) \in S_3$ on $Q_8$ is given by the rule: $(12)
i =-j$, $(12) j =-i$, $(12) k =-k$, and the rules of action of the
other two transpositions are then obvious. In particular, $(123) i
=(23)(12) i= j$, and similarly $(123) j = k$, and $(123) k = i$.
It is then easily checked that $a$ is an involution, that $b$ has
order $8$ and $ab$ has order $3$. Note that this map is given by
identifying antipodal octagons as numbered (and the associated
adjacency of hexagons). Note also that the sixth octagon is
omitted from this picture, but occurs as the outer edges of the
outer hexagons. We show a tree of hexagons in this map, whose
boundary is a cycle missing only two vertices. The left picture
shows this same tree in the corresponding hexagon graph, the
Moebius-Kantor graph of order $16$, and the right picture shows a
Hamilton tree of faces, including also an octagon, whose boundary
is a Hamilton cycle in this map.}
%{\rm In Figure~\ref{fig:q8s3}, we show the
%the Cayley graph of the group $G = Q_8 \rtimes S_3$
%with  a $(2,8,3)$-presentation $\la a,b | a^2=b^8=(ab)^3 =1, etc. \ra$,
%where $a=(12)(37)(56)$ and $b=(12345678)$ live inside $S_8$, the corresponding Cayley map,
%and the Moebius-Kantor graph, the corresponding hexagon graph.}
\end{example}

\bigskip
\bigskip

%12
%\bigskip

%\begin{figure}[h]
%\begin{center}
%\label{fig:q8s3}
%\includegraphics[width=6cm]{Figure4a.eps}
%\caption{cetrta slika}
%\end{center}
%\end{figure}

\bigskip
\bigskip

\begin{figure}[h!]
\begin{center}
%\label{fig:a4}
\includegraphics[width=15cm]{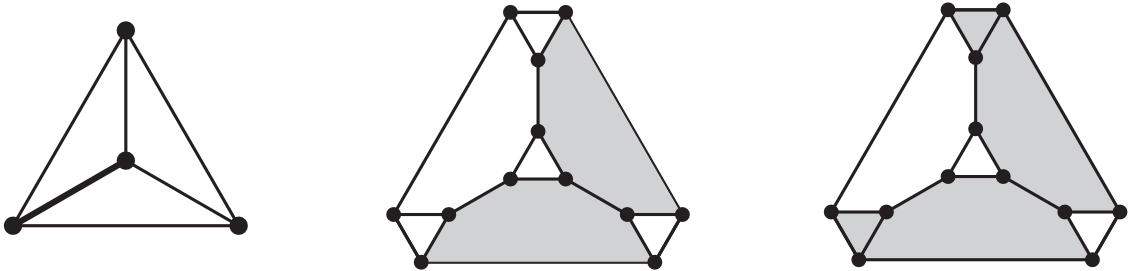}
\caption{A tree of faces in the spherical Cayley map of a Cayley
graph of $A_4$ giving rise to a cycle missing two vertices, the
associated hexagon graph, and a modified Hamilton tree of
faces.}\label{fig:a4}
\end{center}
\end{figure}

\begin{example}
\label{ex:a4}
{\rm In the middle picture of  Figure~\ref{fig:a4} we show
a tree of hexagons, whose boundary is a cycle missing only two vertices
in the spherical Cayley map of
a Cayley graph $X$ of the group $G = A_4$
with  a $(2,3,3)$-presentation $\la a,b \mid a^2=b^3=(ab)^3 =1 \ra$,
where $a=(12)(34)$ and $b=(123)$.
The left picture shows this same tree in the corresponding
hexagon graph $K_4$,
and the right picture shows a Hamilton tree of faces,
including also two triangles, whose boundary is a Hamilton cycle
in this map.}
\end{example}

\bigskip
\bigskip

\bigskip
\bigskip

\begin{figure}[h!]
\begin{center}
%\label{fig:a5}
\includegraphics[width=15cm]{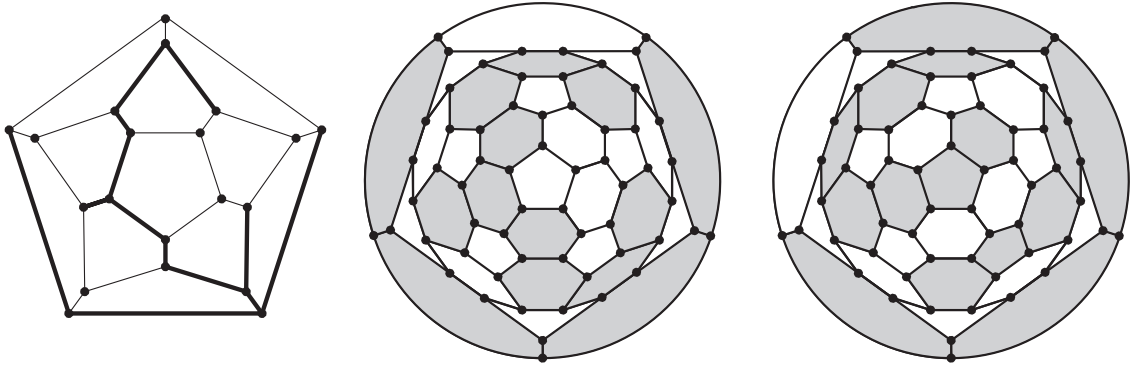}
\caption{A tree of faces in the spherical Cayley map of a Cayley
graph of $A_5$ giving rise to a cycle missing two vertices, the
associated hexagon graph, and a modified Hamilton tree of
faces.}\label{fig:a5}
\end{center}
\end{figure}

\begin{example}
\label{ex:a5}
{\rm In the middle picture of  Figure~\ref{fig:a5} we show
a tree of hexagons, whose boundary is a cycle missing only two vertices
in the spherical Cayley map of
a Cayley graph $X$ of the group $G = A_5$
with  a $(2,5,3)$-presentation $\la a,b \mid a^2=b^5=(ab)^3 =1 \ra$,
where $a=(12)(34)$ and $b=(12345)$.
The left picture shows this same tree in the corresponding
hexagon graph, the dodecahedron,
and the right picture shows a Hamilton tree of faces,
including also two pentagons, whose boundary is a Hamilton cycle
in this map.}
\end{example}

\bigskip
\bigskip
%[h,t]

\bigskip

%
%
%%%%%%%%%%%%%%%%%%%%%%%%%%%%%%% section 3 %%%%%%%%%%%%%%%%%%%%%%%%%%%%%%%%%%%%%
%
%

\section{Cyclic stability and cyclic connectivity}
\label{sec:cyclic}
\noindent

\medskip

A successful application of the method described in the previous section
depends heavily on two purely graph-theoretic results.
The first one, due to Payan and Sakarovitch \cite{PS75}, goes back to 1975 and
deals with maximum sizes of vertex subsets in cubic graphs
inducing acyclic subgraphs, whereas the second one,
due to Nedela and \v Skoviera \cite{NS95}, is somewhat more recent
and concerns cyclic connectivity of vertex-transitive graphs.

Following \cite{PS75}, a paper that is presumably not readily available,
we say that, given a graph (or more generally a loopless multigraph)
$X$, a subset $S$ of $V(X)$ is {\em cyclically stable} if the induced subgraph $X[S]$
is acyclic (a forest). The size $|S|$
of a maximum cyclically stable subset $S$ of $V(X)$
is said to be the {\em cyclic stability number} of $X$.
The following result giving an upper bound
on the cyclic stability number is due to Jaeger \cite{FJ74}.
For the sake of completeness we include its proof.

\begin{proposition}
\label{pro:upper}
{\rm [Jaeger, 1974]}
Let $X$ be a cubic loopless multigraph of order $n$ and let
$S$ be a maximum cyclically stable subset of $V(X)$.
Then

\begin{equation}
\label{equ:cycstab}
|S| = (3n-2c-2e)/4,
\end{equation}

\noindent
where $c$ is the number of connected components (trees)
in $X[S]$ and $e$ is the number of edges in $X[V(X)\setminus S]$.
In particular, $|S| \leq (3n-2)/4$.
\end{proposition}

\proof
Let $V =V(X)$.
First, in view of maximality of $S$ we
have that a vertex in $V\setminus S$
has at most one neighbor in $V \setminus S$,
so that each of the $e$ edges in $X[V\setminus S]$
is an isolated edge. Now let $f$ and $g$ denote,
respectively, the number of edges in $X[S]$ and the number of edges
with one endvertex in $S$ and the other in $V \setminus S$.
Then we have that $f = |S|-c$ and $g = |S| +2c$.
Of course, $e+f+g = 3n/2$ and so $e+2 |S| +c = 3n/2$,
giving us the desired expression for $|S|$.
Now clearly, the maximum value for $|S|$ occurs when
$e=0$, that is when $V \setminus S$ is an independent set of vertices,
and when at the same time $c=1$, that is when $X[S]$ is a tree.
\Qed

\bigskip

In order to explain the result of
Payan and Sakarovitch, we need to introduce
the concept of cyclic connectivity.
Let $X$ be a connected graph.
A subset $F \subseteq E(X)$ of edges of $X$
is said to be {\em cycle-separating} if $X-F$  is disconnected
and at least two of its components contain cycles.
We say that $X$ is {\em cyclically $k$-connected}, in short {\em c.k.c.},
if no set of fewer than $k$ edges is cycle-separating in $X$.
Furthermore, the {\em edge cyclic connectivity}
$\zeta(X)$ of $X$ is the largest integer $k$ not exceeding
the Betti number $|E(X)| -|V(X)| +1$
of $X$ for which $X$ is cyclically $k$-edge connected.
(This distinction is indeed necessary as, for
example the theta graph $\Theta_2$, $K_4$ and $K_{3,3}$ possess
no cycle-separating sets of edges and are thus
cyclically $k$-edge connected for all $k$, however their
edge cyclic connectivities are $2$, $3$ and $4$, respectively.)

In \cite[Th\'eor\`eme~5]{PS75}, Payan and Sakarovitch proved that
in a cubic cyclically $4$-connected graph the above
upper bound for its cyclic stability number given in Proposition~\ref{pro:upper}
is always attained. More precisely, bearing in mind the expression for
the cyclic stability number given in
(\ref{equ:cycstab}), the following result may be deduced from
\cite[Theoreme~5]{PS75}.

\begin{proposition}
\label{pro:pasa}
{\rm [Payan, Sakarovitch, 1975]}
Let $X$ be a cyclically $4$-connected cubic graph of order $n$,
and let $S$ be a maximum cyclically stable subset of $V(X)$.
Then $|S| = [(3n-2)/2]$ and more precisely, the following hold.
\begin{itemize}
\item[(i)]  If $n \equiv 2\,(mod \,4)$ then $|S| = (3n-2)/4$,
and $X[S]$ is a tree and $V(X)\setminus S$ is an independent set of vertices;
\item [(ii)] If $n \equiv 0\,(mod \,4)$ then $|S| =(3n-4)/4$,
and either $X[S]$ is a tree and $V(X)\setminus S$ induces a graph with a single edge,
or $X[S]$ has two components and $V(X)\setminus S$ is an independent set of vertices.
\end{itemize}
\end{proposition}

The connection between cyclic stability and hamiltonicity
is now becoming more transparent.
Let $G$ be a group with a $(2,s,3)$-presentation, $X$ be
the corresponding Cayley graph
and $Y=Hex(X)$ be its hexagon graph .
As described in the previous section it is precisely
the fact that one is able to decompose the vertex set
$V(Y)$ into two subsets, the first one
inducing a tree, and its complement being an independent set of vertices
that enabled us to produce a Hamilton cycle in the original
graph $X$ for the $2,6,3)$-presentations of $\ZZ_6$ and $S_3 \times \ZZ_3$.
Further, with a slight modification,
when the decomposition is such that the first set induces
a tree and its complement induces a subgraph with a single edge,
then a long cycle missing  only two vertices is produced
in $X$. Therefore if $|G|$,
and hence the order of the hexagon graph $Hex(X)$,
is congruent to $2$ modulo $4$,
then part (i) of Proposition~\ref{pro:pasa}
does the trick, provided of course that the $Hex(X)$
is indeed a $c.4.c.$ graph.
On the other hand,
if the $|G|$,
and hence the order of $Hex(X)$,
is divisible by $4$,
then we are not quite there yet
for only one of the possibilities given
in part (ii) of Proposition~\ref{pro:pasa}
will allow us to construct a long cycle in the original graph $X$.
In what follows we explore this situation by, first,  bringing into
the picture an important result on cyclic edge
connectivity of cubic graphs due to Nedela and \v Skoviera,
and second, by showing that, save for a few exceptions,
a cyclically stable set in a hexagon graph of order divisible by $4$,
may always be chosen in such a way that it induces a tree, and its complement
induces a subgraph with a single edge.

The following result is proved in \cite[Theorem~17]{NS95}.

\begin{proposition}
\label{pro:nesko}
{\rm [Nedela, \v Skoviera, 1995]} The cyclic edge connectivity $\zeta(X)$
of a cubic vertex-transitive graph $X$ equals its girth $g(X)$.
\end{proposition}

Consequently, the cyclic edge connectivity
of the hexagon graph, which is an arc-transitive, and thus
also a vertex-transitive cubic graph, coincides with its girth.
As it will soon become clear, a lower bound on the girth is needed.
This is what the next proposition does, where we show that,
with a few exceptions,
the girth of such a graph is not less than $6$.

\begin{proposition}
\label{pro:girth}
Let $X$ be a cubic arc-transitive graph. Then one of the following occurs.
\begin{itemize}
\item [(i)]The girth $g(X)$ of $X$ is at least $6$; or
\item [(ii)] $X$ is one of the following graphs: the theta graph $\Theta_2$, $K_4$,
$K_{3,3}$, the cube $Q_3$, the Petersen graph $GP(5,2)$
or the dodecahedron graph $GP(10,2)$.
\end{itemize}
\end{proposition}

\proof
Clearly, $\Theta_2$ is the only arc-transitive cubic (multi)graph
of girth 2.

Let $G=\Aut X$.
Suppose first that $g(X)=3$. Let $v \in V(X)$ and let $u_0$,
$u_1$ and $u_2$ be its neighbors. By arc-transitivity there exists
an automorphism $\alpha$ of $X$ fixing $v$ and cyclically permuting its neighbors,
that is, $\alpha(u_i) = u_{i+1}$, $i \in \ZZ_3$.
Since $g(X)=3$ it clearly follows that each $u_i$ is adjacent to the other two
neighbors of $v$, and so $X \cong K_4$.

Suppose next that $g(X)=4$. Let
$v \in V(X)$, $N(v) = \{u_i | i \in \ZZ_3\}$
and $\alpha \in G_v$ have the same meaning as above.
Since $g(X) =4$, there are no edges in $N(v)$,
but there must exist, say,  a vertex $x_{01}$,
which is adjacent to both $u_0$ and $u_1$.
If $x_{01}$ is also a neighbor of $u_2$, then it is easily seen
that there exists a third common neighbor of $u_0$, $u_1$ and $u_2$,
implying that $X \cong K_{3,3}$.
If on the other hand, $x_{01}$ is not adjacent to $u_2$,
then there must exist vertices $x_{12}$ and $x_{20}$
which are common neighbors of, respectively,
$u_1$ and  $u_2$,
and of $u_2$ and $u_0$.
But then, using the fact that $X$ is an arc-transitive graph of girth $4$,
one can easily show that the three vertices
$x_{01}$, $x_{12}$ and $x_{20}$ have a common neighbor,
forcing $X$ to be the cube $Q_3$.

Finally, suppose that $g(X)=5$.
We show first that the order of $G$ is divisible by $5$.
Let $v \in V(X)$, let $N(v) = \{u_i | i \in \ZZ_3\}$
be its neighbors' set. We may assume that the 2-arc $u_0vu_1$
is contained on a $5$-cycle.
By arc-transitivity the arc $vu_2$ must also be
contained on a 5-cycle,
and so either the 2-arc $u_0vu_2$ or the 2-arc $u_1vu_2$
is contained on a 5-cycle. But then both 2-arcs are on a 5-cycle.
In short, an edge of $X$ is contained on at least two distinct $5$-cycles.
On the other hand, an edge of $X$ cannot be contained on more than four
$5$-cycles. Namely, since $g(X)=5$,
a $3$-arc of $X$ is contained on at most one $5$-cycle,
and so a $2$-arc of $X$ is contained on at most two $5$-cycles,
and a $1$-arc on at most four $5$-cycles.
More precisely, either each $3$-arc in $X$ gives rise
to a unique $5$-cycle or each $2$-arc gives rise to a unique $5$-cycle.
In the first case there is a
total of $6n/5$ cycles of length $5$ in $X$
(with each edge on four $5$-cycles),
and in the second case there is a total of $3n/5$
cycles of length $5$ in $X$
(with each edge on two $5$-cycles), where $n$ is the order of $X$.
In both cases  we have that $5$ divides $n$.
In particular, $X$ has an automorphism of order $5$.

Next we show that the stabilizer of the action of $G$ on the set $\C$ of all
$5$-cycles in $X$ contains an element of order $5$.
In other words, we show that there is an automorphism
of $X$ which rotates a $5$-cycle in $X$ (and so
in the terminology of \cite{NB78},
$X$ has a consistent $5$-cycle).

In view of the above remarks on the numbers of $5$-cycles
an edge of $X$ is contained on, we can easily see
that either $G$ is transitive on $\C$
or $G$ has two equal length orbits, say $\C_1$ and $\C_2$
in its action on $\C$. But note that  $|\C|$ is either
$3n/5$ or $6n/5$ and moreover $|G| =3\cdot2^r n$,
by the classical result of Tutte \cite{Tut47}
(where $r\leq2$, as $G$ is at most $3$-arc-transitive).
Let $P$ be a Sylow $5$-subgroup of $G$, say of order  $5^k$.
It follows that there must be an orbit in the action of
$P$ on $\C$ (or $\C_1$ and $\C_2$) of length
$5^t\cdot m$, where $(5,m)=1$ and $t <k$.
This implies that there is a stabilizer in the action of $G$
on $\C$ (or $\C_1$ and $\C_2$) containing an element of $P$.
In other words, there is an automorphism of $X$ rotating a $5$-cycle.
This will prove crucial in the final steps of our proof.

Now let $\rho$ be this automorphism of order $5$ in $X$
rotating a $5$-cycle $C = v_0v_1v_2v_3v_4$ of $X$.
Clearly, since $X$ is cubic the action of $\rho$ on $V(X)$ is semiregular.
Let $u_i$, $i \in \ZZ_5$, respectively, be the additional
neighbors of $v_i$, $i \in \ZZ_5$ so that
$\rho(u_i) = u_{i+1}$ for all $i \in \ZZ_5$.
Supposing first that $U=\{u_i | i \in \ZZ_5\}$ induces
a cycle, one may easily deduce that $X$ is isomorphic to the Petersen graph $GP(5,2)$.
We may therefore assume that $U$ is an independent set of vertices.
Now if the additional neighbors of $u_i$, $i \in \ZZ_5$, were
in two orbits of $\rho$, then no edge $v_iu_i$ would be contained on
a $5$-cycle. Hence there is a third orbit $W=\{w_i | i \in \ZZ_5\}$
of $\rho$, with $\rho(w_i) = w_{i+1}$ for all $i\in\ZZ_5$,
containing all of the additional neighbors of vertices in $U$.
Of course, $W$ is an independent set of vertices,
and there is a fourth orbit $Q = \{q_i | i \in \ZZ_5\}$
of $\rho$ containing the additional neighbors of vertices in $W$.

But now in order for the edges with one endvertex in $W$
and the other endvertex in $Q$ to lie on a $5$-cycle,
the orbit $Q$ must necessarily induce a cycle.
In other words, $\rho$ has precisely four orbits
and so $X$ is a non-bipartite cubic arc-transitive graph of order $20$.
Hence $X$ is isomorphic to the graph of the dodecahedron,
that is the graph $GP(10,2)$ in the generalized Petersen graph notation.
\Qed

\bigskip

Using the previous two results we will
give a refinement to Payan and Sakarovitch result
in Proposition~\ref{pro:pasa}
by showing that for cubic arc-transitive graphs
of girth $6$ the maximum cyclically stable subset
may always be chosen to induce a tree.

\begin{proposition}
\label{pro:modi}
Let $X$ be a cubic arc-transitive graph of order $n \equiv 0\, (mod \, 4)$,
not isomorphic to any of the following graphs:
$K_4$, the cube $Q_3$, or the dodecahedron graph $GP(10,2)$.
Then there exists a cyclically stable subset
$S$ of $V(X)$ which induces a tree, and such that $V(X) \setminus S$
induces a graph with a single edge.
\end{proposition}

\proof
Observe that, in view of Proposition~\ref{pro:girth},
the girth of $X$ is at least $6$, and hence,
in view of Proposition~\ref{pro:nesko},
$X$ is  a $c.6.c.$ graph.
Note also that the statement of this proposition
really says that in part (ii) of Proposition~\ref{pro:pasa},
a particular one of the two possibilities may be chosen.

We procceed as follows.
Let us first modify our graph $X$ by deleting a pair of adjacent vertices,
say $u$ and $v$. This modified graph $Y = X -\{u,v\}$ has
$n-2$ vertices with the two neighbors $u_1$ and $u_2$ of $u$ and
the two neighbors $v_1$ and $v_2$ of $v$
having valency $2$, and all the remaining vertices having
valency $3$. ``Forgetting'' the four vertices $u_1$, $u_2$, $v_1$ and $v_2$,
we may therefore also think of this
modified graph as being cubic of order $n-6$.
The important thing however is that $Y$ must be a $c.4.c.$ graph.
Namely, taking two vertex disjoint
cycles $C_1$ and $C_2$ in
$Y$, a maximum number of additional paths separating these two cycles
one can obtain by adding the vertices $u$ and $v$
(and all their neighbors), is $2$.
This would occur if each of the two vertices $u$ and $v$ had one
of their two (additional) neighbors in $C_1$ and the other in $C_2$.
We conclude that by going from $Y$ back to the original graph $X$,
the cyclic connectivity can go up by at most $2$.
Since $X$  is $c.6.c$ it follows that $Y$ is $c.4.c$.

But the order of $Y$ is congruent to $2$ modulo $4$
and by part (i) of Proposition~\ref{pro:pasa}
(in this particular instance we are forgetting the four
neighbors to make $Y$ cubic!!),
there exists a maximum cyclically stable subset $R$
of $V(Y)$ inducing a tree and
such that its complement $V(Y) \setminus R$ is an independent set
of vertices. Now the maximum cyclically stable subset $S$ of $V(X)$
is now obtained by taking $S = R \cup \{u_1,u_2,v_1,v_2\}$.
(The edge $uv$ is thus the single edge of the graph induced on the complement
$V(X) \setminus S$.)
\Qed

%
%
%%%%%%%%%%%%%%%%%%%%%%%%%%%%%%% section 4 %%%%%%%%%%%%%%%%%%%%%%%%%%%%%%%%%%%%%
%
%

\section{Proving Theorem~1.1}
\label{sec:proof}
\noindent

\medskip

\bigskip
{\bf Proof of Theorem~\ref{the:main}.}
Suppose first that $|G| \equiv 2\,(mod\,4)$.
As demonstrated in Section~2
a tree in $Hex(X)$ whose complement is an independent set
of vertices gives rise to a Hamilton tree of faces
in the Cayley map associated with $X$ and thus to
a Hamilton cycle in $X$. In view of
Example~\ref{ex:z6}, which takes care of the case $Hex(X) \cong \Theta_2$,
and Propositions~\ref{pro:pasa}, \ref{pro:nesko}
and \ref{pro:girth} which combined together take care of
the case $Hex(X) \not\cong \Theta_2$, the graph $X$ is then clearly hamiltonian.
As for the case when $|G| \equiv 0\,(mod\,4)$,
we use Examples~\ref{ex:a4}, {\ref{ex:s4} and \ref{ex:a5}
for the case when $Hex(X)$ is isomorphic, respectively,
to one of $K_4$, $Q_3$ or $GP(10,2)$
and Proposition~\ref{pro:modi} for the case when
$Hex(X) \not\cong K_4, Q_3, GP(10,2)$, to ensure
the existence in the hexagon graph $Hex(X)$ of a
a tree whose vertex complement (in $Hex(X)$) is a graph with  a single edge.
This tree then translates in $X$ into a tree of faces in the Cayley map
(associated with $X$) whose boundary misses only two vertices.
Namely, those two adjacent vertices in $X$
with the corresponding edge being shared
by the two hexagons which (in $Hex(X)$)
are the endvertices of the single edge in the complement
of the three chosen above. Consequently,
$X$ contains a long cycle missing only two vertices.
In particular, $X$ has a Hamilton path in this case.
This completes the proof of Theorem~\ref{the:main}.
\Qed

\bigskip
\noindent
{\bf\large Acknowledgment.}
The authors wish to thank Yuqing Chen, Marston Conder, Gorazd Lah,
Caiheng Li, Klavdija Kutnar, Denis Sjerve, Ron Solomon, and
Tsuyi Yang for helpful conversations about the material in this paper.

\end{document}